\documentclass[12pt]{amsart}
\newfont{\sheaf}{eusm10 scaled\magstep1}
\usepackage{amssymb} 
\usepackage[mathscr]{eucal} 
\usepackage[all]{xy}
\usepackage{amsmath} 

\newcommand{\C}{\ensuremath{\mathbb{C}}}
\newcommand{\R}{\ensuremath{\mathbb{R}}}
\newcommand{\Z}{\ensuremath{\mathbb{Z}}}

\newcommand{\hol}{\ensuremath{\mathcal{O}}}

\newcommand{\PP}{\ensuremath{\mathbb{P}}}

\newcommand{\HHH}{\ensuremath{\mathcal{H}}}

\newcommand{\Proof}{{\it Proof. }}
\newcommand{\ra}{\ensuremath{\rightarrow}}

\newcommand{\sT}{{\mathcal T}}

\newcommand{\la}{\lambda}

\newcommand{\Ga}{\Gamma}

\newcommand{\La}{\Lambda}
\newcommand{\ga}{\gamma}

\newcommand{\QED}{\hspace*{\fill}$Q.E.D.$}

\def\eea{\end{eqnarray*}}
\def\bea{\begin{eqnarray*}}

\def\hol{{\mathcal{O}}}

\newtheorem{teo}{Theorem}[section]
\newtheorem{df}[teo]{Definition}
\newtheorem{lem}[teo]{Lemma}

\newtheorem{oss}[teo]{Remark}
\newtheorem{prop}[teo]{Proposition}

\newtheorem{REM}[teo]{Remark}

\def\eea{\end{eqnarray*}}
\def\bea{\begin{eqnarray*}}

\begin{document}

\title{ Real structures on torus bundles and their deformations}

\author{Fabrizio Catanese\\
     Paola Frediani\\
     }

\footnote{
The present research took place in the framework of the Schwerpunkt
"Globale Methoden in der komplexen Geometrie", and of the PRIN 2003:
"Spazi dei moduli e
     teoria di Lie" (MURST)

 }
\maketitle

\begin{abstract}
We describe the family of real structures $\sigma$ on principal holomorphic torus bundles
$X$ over  tori, and prove its connectedness when the complex dimension is at most 
three.
Hence follows  that  the differentiable type (more precisely, the
orbifold fundamental group) determines the deformation type  
of the pair $(X, \sigma)$ provided we have complex
dimension at most three, fibre dimension one, and a certain 'reality'  condition on the
fundamental group is satisfied.
\end{abstract}

\section{Introduction}
Given an oriented compact differentiable manifold $M$, a quite  general problem
is to consider the space $\mathcal C (M)$ of compatible complex structures
$J$ on $M$.  Its connected components are called "Deformation classes in the large
(of complex structures on $M$)".

A complete answer is known classically for curves, and, by the work of Kodaira,
extending Enriques' classification,
for special algebraic surfaces.

Since $\mathcal C (M)$ is infinite dimensional, one can equivalently 
consider instead the finite
dimensional complex analytic space $ \frak T (M)$, called 
Teichm\"uller space, corresponding
to the quotient of  $\mathcal C (M)$ by the group $ Diff ^0(M)$ of 
diffeomorphisms
isotopic to the identity. Its local structure is described by the 
Kuranishi theory (\cite{ku2}).

Inside $ \frak T (M)$  one has the open set
$\frak KT (M)$ of complex structures admitting a K\"ahler metric.

Hironaka showed (\cite{hir62}) that  $\frak KT (M)$ fails to be a 
union of connected components
of  $ \frak T (M)$.

Later, Sommese showed (\cite{somm75}) that, even in the case of
complex tori  of dimension
$\geq 3$,  $ \frak T (M)$ has more than
one connected component. It was known to Kodaira and Spencer that in this case
$\frak KT (M)$ is connected and coincides with the subset of 
translation invariant complex
structures.

It seems to be  more the rule than the exception that   $ \frak T (M)$ has
a lot of connected components. For instance, already for surfaces of 
general type
it was shown by Manetti (\cite{man01}), and by Catanese and Wajnryb 
in the simply connected
case (\cite{c-w04}), that the number of connected components of  $ 
\frak T (M)$ may be arbitrarily
large.

The above  work shows also that perhaps a more approachable 
version of the general problem is the following:

QUESTION  :  given a complex structure 
$ X : = (M, J)$,  determine the connected component $\mathcal C^0  (X)$ of 
$\mathcal C (M)$ containing it.  I.e., given a complex manifold $X$, 
find its deformations in  the large.

In \cite{cat02} (cf. also \cite{cat04})
it was shown that every deformation in the large of a complex
torus is a complex torus.  This was used to determine the stability 
by deformations in the large of
holomorphic torus bundles over curves of genus at least two ( \cite{cat04}),
and to obtain the same result  (\cite{cf}) for threefolds which 
are  holomorphic torus bundles
over  two dimensional tori, provided that the
fundamental group satisfies a suitable 'reality' 
condition.

Consider now, as customary, a real manifold as a pair $(X, \sigma)$, where $X$ is a 
compact complex manifold
and $\sigma$ is a antiholomorphic involution.

\begin{df}
$(X, \sigma)$ is said to be {\bf simple  } if, once  we fix the 
differentiable type of the pair $(X, \sigma)$, then we have a unique deformation class.

I.e., $M$ being the oriented differentiable manifold underlying  $X$,
simplicity holds iff  the space $\mathcal R (M,\sigma)$ of compatible
complex structures (which make $\sigma$ antiholomorphic) is
connected.

$(X, \sigma)$ is said to be {\bf quasi- simple  } if the above holds when we
restrict ourselves to complex structures which  are deformations of $X$
(cf. \cite{dk}). I.e., $(X, \sigma)$ is quasi- simple   iff  $\mathcal R (M,\sigma)
\cap \mathcal C ^0  (X)$ is connected.

More generally, as pointed out to us by Itenberg (cf. \cite{dik}), one can pose the problem of
 quasi-simplicity for the more general  case 
of a finite subgroup $G$ of the Klein group of dianalytic
(i.e., biholomorphic or antibiholomorphic)
automorphisms of $X$ ( $ G \cong \Z / 2$ in the previous case).

\end{df}

A basic question here is : which real manifolds are simple ? Which ones are quasi simple ?

The answer is positive for curves, by work of Sepp\"ala and Silhol (\cite{s-s}), 
and  also for many
surfaces of special type, thanks to work of many authors, like 
Comessatti, Silhol,
Nikulin, Degtyarev, Itenberg, Kharlamov,  Catanese-Frediani, Welschinger.
These works (cf. references) give evidence to the following

CONJECTURE 3.2  Special surfaces are simple.

Motivated by this conjecture,
in \cite{c-f03} we introduced the notion of  orbifold fundamental group 
in order to prove the
simplicity of  real hyperelliptic surfaces and in order to  achieve a 
complete classification of them
( there are 78 types).

The orbifold fundamental group of a real variety $(X,\sigma)$ is
an exact sequence

$$ 1 \ra \pi_1 (X) \ra  \Pi_{\sigma}: = \pi_1^{orb} (X) \ra \Z/2 \ra 1,$$
where, if the action of $\sigma$ is fixpoint free, 
$\Pi_{\sigma}$ is the fundamental group of the quotient $ X / \sigma$.

We are mainly interested here in the special case where $X$ is 
a $K(\pi , 1)$ and more precisely where the universal cover  $\tilde{X}$ of $X$ 
is contractible. In this case we have
$$ X =  \tilde{X} / \pi_1 (X) ,  \  X / \sigma  = \tilde{X} / \Pi_{\sigma} \ .$$

In the even more special case where $ \tilde{X}$ has no moduli, the simplicity or quasi
simplicity question is translated into the question whether  the conjugacy class of the
embedding 
$ \Pi_{\sigma} \ra Dian ( \tilde{X})$ is parametrized by a connected
variety.

Now, the cases which we treated already, real hyperelliptic surfaces and
Kodaira surfaces (\cite{c-f03} , \cite{f}) are both torus bundles over tori,
whose study is reduced (\cite{cat04}) to the study  of P.H.T.B.T. = principal holomorphic
torus bundles over tori. The basic question which is of interest to us
is  to determine when a real P.H.T.B.T. is simple.

One could consider in greater generality the following situation:
we have a nilpotent group $\Pi$, such that the factors of the central series
are torsion free abelian of even rank.

Then the  real Lie group  $\Pi \otimes \R$ contains $\Pi$ as
a cocompact lattice (i.e., as a discrete subgroup with compact quotient $ M : = 
(\Pi \otimes \R) / \Pi$) and we consider the differentiable manifold
$ M : =  (\Pi \otimes \R )/ \Pi$. 

Already the case of tori has warned us that there could be too many
exotic complex structures, so we restrict ourselves to consider
only the right invariant complex structures on the homogeneous space $M$.

Obviously, the  right invariant  almost complex structures $J$ correspond
exactly
to the complex structures on the tangent space at the origin, but
there remains to write down explicitly the integrability conditions
for a given almost complex structure $J$.

This is rather simple in the case of  P.H.T.B.T. , treated in \cite{cat04}
and in \cite{cf}: we have an analogue of the Riemann bilinear relations
(classically known for line bundles on complex tori)  which are  
easy to explain as follows (cf.  the next section for more details).

Consider the central extension of fundamental groups given
by the homotopy exact sequence of the fibre bundle  $f : X \ra Y$:
$$  1 \ra \pi_1(T):  = \La \ra \Pi := \pi_1 (X) \ra \pi_1(Y) := \Ga \ra 1.$$

The central extension is classified by an element
$$ A \in H^2 ( \Ga, \La ) \cong \La^2 (\Ga^{\vee} )\otimes \La .$$

As differentiable manifolds the base $ Y = (\Ga \otimes \R  )/ \Ga$, and the fibre
$  T=( \La \otimes \R ) / \La$,  acquire a varying complex structure
by writing a Hodge decomposition
$$ \Ga \otimes \C = V \oplus \bar{V},  \ \ \La \otimes \C = U \oplus \bar{U}.$$

The Riemann bilinear relations (which are equivalent to the integrability of the given
almost complex structure on $ M : =( \Pi \otimes \R )/ \Pi$) amount
to saying that the component of $A$ in $\La^2 (V^{\vee} )\otimes  \bar{U}$
is zero.
These equations characterize then the so-called Appell Humbert 
family of such
bundles, which parametrizes all the right invariant complex structures
on $ M : =( \Pi \otimes \R )/ \Pi$.

The Riemann relations allow a decomposition 
$$  A =  B + \overline{B}, \ \ B = B' + B'',  \  B' \in \Lambda^2(V)^{\vee}
\otimes U, \ B'' \in (V^{\vee} \otimes \overline{V})^{\vee} \otimes U.$$ 

$B'$ is called the holomorphic part, $ B''$ is called the Hermitian part, and
this decomposition  reduces to
questions of multilinear algebra the description of the spaces 
of holomorphic forms on $X$, of the subspace of closed holomorphic forms.
 and the explanation of an important phenomenon discovered by Nakamura
(\cite{nak75}), namely the fact that local deformations of a parallelizable manifold need 
not be parallelizable.

One has for instance the following

 {\bf Theorem, \cite{cat04}:} $X$ is parallelizable if and only if $ B'' = 0$.

The Riemann  bilinear relations define a parameter space which parametrizes 
the complex structures in the Appell Humbert family, and which is denoted
by  ${\mathcal T}'{\mathcal B}_{A} $, since it parametrizes torus bundles.

The parameter space ${\mathcal T}'{\mathcal B}_{A} $ is smooth if
$ m = dim_{\C} (Y) \leq2, d = dim_{\C} (T)  = 1$ (and it is singular already
for $d=1, m \geq 3$, at the points where $ B '' = 0$).  

The smoothness
of  ${\mathcal T}'{\mathcal B}_{A} $ makes it possible to analyse the versality of the
(complete) Appell-Humbert family  and to prove it (\cite{cf}) under the
following\\
 
{\bf  'Reality' condition for $ A $ in the case  $d=1$ } :
choose a basis $ e_1, e_2$ of $\La$ and write $ A = A_1 \otimes e_1 +
A_2 \otimes e_2$. Then we want a  real solution for the equation
$$ Pfaffian  (\la_1 A_1 + \la_2 A_2 ) = 0  . $$ 

{\bf Theorem B , \cite{cf}}{\em  The Appell Humbert family yields an open
subset of  $\mathcal C ( (\Pi \otimes \R)/ \Pi)$ if $ d=1, m=2$ and moreover

i) $A$ is non degenerate

ii)  $ dim ( Im A) = 2 $

iii) $A$ satisfies the 'reality' condition.}

Combined with 

{\bf Theorem A , \cite{cat04}} {\em  The Appell Humbert family yields a closed
subset of  $\mathcal C ( (\Pi \otimes \R)/ \Pi)$ if $ d=1$ and moreover
 $ dim ( Im A) = 2 $.}

we find that, under the hypotheses of Theorem B, the Appell Humbert family
yields a connected component of  $\mathcal C ( (\Pi \otimes \R)/ \Pi)$.

In this paper we describe completely,  via explicit
equations, the family which
parametrizes  real structures on
a torus bundle in the Appell  Humbert family 
 and then we show the connectivity of the family in special cases.

We obtain for instance  the following result, which we formulate here
in the less technical form

{\bf Theorem C }{\em   Assume again $ d=1, m=2$, that $A$ is nondegenerate, and 
fix the orbifold fundamental group $\Pi_{\sigma}$.
Then the real structures for torus bundles in the Appell  Humbert family
are parametrized by a connected family.}

We obtain therefore as a corollary

{\bf Theorem D }{\em   Same assumptions as in theorem B: $ d=1, m=2$, and 
 i), ii), iii) are satisfied.

Then simplicity holds for  real  torus bundles in the Appell  Humbert family.}

The paper is organized as follows : in section 2 we recall the theory developed
in \cite{cat04} and \cite{cf} concerning principal holomorphic torus bundles
over tori, while section 3 is devoted to recalling
the results already established for  the Appell Humbert family.

The bulk of the paper is section 4, which studies the real structures on torus
bundles in the Appell Humbert family.  Strangely enough, the  extra symmetries
coming  from the existence
of a real structure makes calculations  somewhat easier.

In fact, the real structure allows  splittings
$$\Ga \otimes \R = V^+ \oplus V^- \ , \ \La \otimes \R = U^+ \oplus U^- ,$$

and, after we show that the orbifold fundamental group   $\Pi_{\sigma}$
has a dianalytic affine representation in the complex vector space
$( V \otimes \C ) \oplus ( U \otimes \C ) $, we can describe the complex structures
in the Appell Humbert family which  are compatible with the real structure
as a pair of linear maps
$$ B_2 : V^- \ra V^+ \ , \   B_1 : U^- \ra U^+ \ , $$
satisfying explicit second degree matrix equations.  

To illustrate the power of this way of calculating, we give a short new proof of 
the second author's result that simplicity holds for Kodaira surfaces.

Then we proceed to the case of threefolds, obtaining our main results (Theorems C and D).

\section{Principal holomorphic torus bundles over tori: generalities}

Throughout  the paper, our set up will be the
following: we have a holomorphic submersion between compact
complex manifolds

$$ f : X \ra Y  ,$$
such that the base $Y$ is a complex torus, and one fibre $F$  (whence
all the fibres, by
theorem 2.1 of
\cite{cat04}) is also a  complex torus.

We shall denote this general situation by saying that $f$ is
differentiably a torus bundle.

We let $n= dim X$, $m = dim Y$, $d=dim F= n -m$.

In general  ( cf.\cite {fg65})) $f$ is a
holomorphic bundle if and only if all the  fibres are
biholomorphic.

This holds necessarily in the special case $d=1$, because the moduli space
for $1$-dimensional complex tori exists and is isomorphic to $\C$.

Assume now more specifically that we have a holomorphic torus fibre bundle, thus we have
     (cf. \cite{cat04}, pages 271-273) the exact sequence
of holomorphic vector bundles 
$$ \  0 \ra  \Omega^1_Y \ra f_* \Omega^1_X \ra f_* \Omega^1_{X|Y} \ra 0 .$$
We have a principal
holomorphic bundle if moreover $f_* \Omega^1_{X|Y}$ is a  trivial
holomorphic bundle.

\begin{oss}
In general (cf. e.g. \cite{bpv84}) if $T$ is a complex torus, we have an
exact sequence of complex Lie groups
$$ 0 \ra T \ra Aut(T) \ra M \ra 1 $$
where $M$ is discrete.
Taking sheaves of germs of holomorphic maps with source $Y$ we get
$$ 0 \ra \HHH(T)_Y \ra \HHH (Aut(T))_Y \ra M \ra 1 $$
hence the exact sequence
$$ 0 \ra H^1(Y, \HHH (T)_Y) \ra H^1(Y, \HHH (Aut(T)))_Y \ra H^1(Y,
M)  \ra H^2(Y, \HHH (T)_Y).$$
Since holomorphic bundles with base $Y$ and fibre $T$ are
classified by the cohomology group $H^1(Y, \HHH (Aut(T))_Y)$,
$H^1(Y, M)$
determines the discrete obstruction for a holomorphic bundle
     to be a principal holomorphic bundle.
\end{oss}
In view of this, the study of holomorphic torus bundles is reduced to the study of principal
holomorphic torus bundles.

In the case of a principal holomorphic bundle we
write $ \La : = \pi_1(T ) $,  $ \Ga : = \pi_1(Y ) $
and the exact sequence
$$\   \ra H^0(\HHH(T)_Y) \ra H^1(Y, \Lambda ) \ra  H^1 (Y,
\hol_Y^d) \ra H^1(\HHH(T)_Y) \ra ^c \ra H^2(Y, \Lambda ) $$
determines a  cohomology class
$\epsilon \in H^2(Y, \Lambda )$
which classifies the  central extension
$$  1 \ra \pi_1(T) = \La \ra \Pi := \pi_1 (X) \ra \pi_1(Y) = \Ga \ra 1$$
(it is central by the triviality of the monodromy automorphism).

\begin{prop}

\label{unicover}

\cite{cf}
Let $f: X \rightarrow Y$ be a principal holomorphic torus bundle over
a torus as
above.

Then the universal covering of $X$ is isomorphic to ${\C}^{m+d}$ and $X$ is
biholomorphic to a quotient $X \cong {\C}^{m+d}/\Pi$.

\end{prop}

Let us briefly recall  the classical  way to look at the family
$\sT_m$ of complex
tori of complex dimension $=m$. We fix a lattice $\Ga$ of rank $ 2m$, and
we look at the complex ($m$-dimensional) subspaces $ V \subset \Ga \otimes \C$
such that  $ V \oplus \bar{V} = \Ga \otimes \C$: to $V$ corresponds the
complex torus $ Y_V : = (\Ga \otimes \C )/ (\Ga \oplus \bar{V} )$.
We finally select one of the two resulting connected components by requiring
that the complex orientation of $V$ induces on $\Ga \cong p_V (\Ga)$
a fixed orientation.

We consider similarly  the complex tori $T_U  : = (\La \otimes \C )/ (\La \oplus \bar{U} )$
which can occur as fibres of $f$.

Consider now  our principal holomorphic torus bundle $f : X \rightarrow Y$
over a
complex torus $Y_V$ of dimension $m$, and with fibre a complex torus
     $T_U$ of dimension $d$  and let  $\epsilon \in
H^2(Y, \Lambda) = H^2(\Ga, \Lambda)$ be the cohomology class classifying
the central extension

\begin{equation}
\label{central}
1 \rightarrow \Lambda \rightarrow \Pi \rightarrow \Gamma \rightarrow 1.
\end{equation}

We reproduce an important result from \cite{cf}
\begin{lem}
\label{difaction}
It is possible to "tensor" the above exact sequence with $\R$,
obtaining an exact
sequence of Lie Groups
$$  1 \ra \Lambda \otimes \R \ra \Pi \otimes \R \ra \Gamma
\otimes \R \ra 1$$ such that
$\Pi$ is a discrete subgroup of $\Pi \otimes \R$ and such that
$X$ is diffeomorphic to the quotient manifold 
$$  M:=  \Pi \otimes \R /  \Pi .$$
\end{lem}

\Proof
Consider, as usual, the map
$$A: \Ga \times \Ga \ra \La,$$
$$A(\ga,\ga') =  [\hat{\ga}, \hat{\ga'}] =  \hat{\ga}
\hat{\ga'}(\hat{\ga})^{-1}(
\hat{\ga'})^{-1},$$
where $\hat{\gamma}$ and $\hat{\gamma'}$  are  respective liftings to $\Pi$ of
elements $\ga, \ga' \in \Ga$.
We observe that since the extension (\ref{central}) is central,
     the definition of $A$ does
not depend on the choice of the liftings of $\ga$, resp. $\ga'$ to $\Pi$.

As it is well known, $A$ is bilinear and alternating, so $A$ yields a cocycle in
$H^2(\Ga, \La)$ which "classifies" the central extension (\ref{central}).
Let us review how does this more precisely hold.

Assume  that $\{\ga_1, ..., \ga_{2m}\}$ is a basis of $\Ga$ and choose
     fixed  liftings $\hat{\ga_i}$ of $\gamma_i$ in $\Pi$, for each $i =
1, \dots , 2m$.
Then automatically we have determined a canonical way to lift
elements $\ga \in \Ga$ through:
$$ \ga = {\ga}_1^{n_1} \dots  {\ga}_{2m}^{n_{2m}} \mapsto \hat{\ga} : =
{\hat{\ga}_1}^{n_1}... {\hat{\ga}_{2m}}^{n_{2m}}.$$
Hence a canonical way to write the elements of $\Pi$ as products
     $\lambda \hat{\ga}$,
where  $\la \in \La$ and $\hat{\ga}$ is as above.

Since $\forall i,j$, one has

$$\hat{\ga_i} \hat{\ga_j} = A(\ga_i, \ga_j) \hat{\ga_j} \hat{\ga_i},$$

      we have a standard way of computing the products
$(\lambda \hat{\ga}) (\lambda' \hat{\ga}')$ as
$ \la'' (\widehat{(\ga \ga'}))$,  where $\la ''$ will be
computed using $A$.

We can also view $\Pi$ as a group of affine transformations of $(\La
\otimes {\R})
\oplus (\Ga \otimes \R)$. In fact,   $(\La \otimes {\R})
\oplus (\Ga \otimes \R)$ is a real vector space with basis
$\{\lambda_1,..., \la_{2d},
\ga_1, ..., \ga_{2m}\}$ where $\{\lambda_1,..., \la_{2d}\}$ is a basis of
$\La$ and the action of $\Pi$ on $(\La \otimes {\R})
\oplus (\Ga \otimes \R)$
is given as follows:

$\lambda_i$ acts on $(\La \otimes {\R}) \oplus (\Ga \otimes \R)$
sending $(y, x)$ to
$(y + \lambda_i, x)$, while the action of $\hat{\gamma_j}$ is defined using the
multiplication $(\lambda \hat{\gamma})
\mapsto (\lambda \hat{\gamma}) \hat{\gamma_j}$.

More precisely if $y  \in \La \otimes \R$, $x = \sum x_j \gamma_j \in
\Ga \otimes
\R$, $\ga'  = \sum \nu_h \ga_h \in \Ga $, we have
$$(y,x)\hat{\ga'} : =  (y + \phi_{ \ga'}(x), x + \ga'),$$
where
$$\phi_{ \ga'}(x) = \sum_{j \geq h} x_j \nu_h A( \gamma_j,  \gamma_h) =
\sum_{j \geq h} x_j A_{jh} \nu_h = ^tx T^-  \ga',$$
where $T^-$ is the lower triangular part of the matrix $A$, so that
we can write
$A = T^- - ^tT^- $.

Therefore  we can endow
$(\La \otimes {\R}) \oplus (\Ga \otimes \R) =: \Pi \otimes \R$ with a Lie group
structure defined by
$$(y,x) (y', x') = (y + y' + T^-(x,x'), x + x'),$$
and the quotient $(\Pi \otimes \R)/\Pi$  of this Lie group by the
discrete subgroup $\Pi$ is immediately seen to be
diffeomorphic to
$X$.

     \QED

\begin{REM}
\label{difaction1}
We can change coordinates in $(\La \otimes {\R}) \oplus (\Ga \otimes
\R)$ in such a
way that the action of  the set $\hat{\Ga} \cong \Ga \subset \Ga
\otimes \R$  on $\Pi
\otimes
\R$ is given by
$$(y,x) \hat {\ga} = (y + A(x, \ga) + 2 S ({\ga}, \ga) , x + \ga),$$
where $S ({\ga}, \ga') $ is a symmetric bilinear  $(\frac{1}{4} \La )$-
valued form, and $ 2 S ({\ga}, \ga) \in \La$.
\end{REM}

\Proof

Let us define the symmetric form $S :  = - \frac{T ^-+^tT^-}{4}$, so
that $T^-  + 2S=
\frac{T^- - ^tT^-}{2} = \frac{A}{2}$.

Consider the map $ \Psi : (\La \otimes {\R}) \oplus (\Ga \otimes \R) \ra
(\La \otimes {\R})
\oplus (\Ga \otimes \R)$ defined by $\Psi (y,x) = (2(y + S(x,x)), x) =:( \eta, x)$.

Then $\forall \hat {\ga}  \in \hat {\Ga} $ we have an induced action
$$(\eta, x) \hat {\ga}  = \Psi ((y,x) \hat {\ga} ) = \Psi (y + T^-(x, \ga),
x + \ga) =$$
$$ (2y + 2 T^-(x, \ga) + 2S(x+\ga,x + \ga), x + \ga) =
(\eta + A(x, \ga) + 2S(\ga,\ga), x + \ga),$$
and we conclude observing that $2S(\ga,\ga)   = T^-(\ga,\ga) \in \La$.

\QED

We recall from \cite{cat04}  the First Riemann
bilinear Relation: it is derived from the exact cohomology
sequence
$$ H^1 (Y, \hol_Y\otimes U) \cong H^1 (Y, \HHH(U)_Y) \ra H^1(\HHH(T)_Y)
\ra ^c \ra H^2(Y, \Lambda ) \ra H^2 (Y, \HHH(U)_Y) $$
     and says that the class $\epsilon$ maps to zero in $H^2 (Y, \HHH(U)_Y) $.
More concretely  the

{\underline{\bf First Riemann Relation for  principal holomorphic
Torus Bundles} }

is expressed as follows:

{\bf Let $A : \Gamma \times \Gamma \ra \Lambda$ be the alternating bilinear map
representing the cohomology class $\epsilon$ :
     then
     $$ A \in \Lambda^2(\Gamma
\otimes \R )^{\vee} \otimes (\Lambda \otimes \R ) \subset  \Lambda^2(\Gamma
\otimes \C )^{\vee} \otimes (\Lambda \otimes \C) \subset \Lambda^2(V \oplus
\bar{V})^{\vee} \otimes (U \oplus \bar{U}),$$
satisfies the property that its  component  in $
\Lambda^2(
\bar{V})^{\vee} \otimes (U )$ is zero.}

It is important, for the forthcoming calculations, to understand  in detail the bilinear algebra
underlying the Riemann bilinear
relation.

We observe preliminarly that  one has a natural isomorphism
$ \Lambda^2(V \oplus
\overline{V})^{\vee}  \cong  \Lambda^2(V)^{\vee}  \oplus
(V^{\vee} \otimes \overline{V})^{\vee}   \oplus
\Lambda^2(\overline{V})^{\vee}$,
where the middle summand embeds by the wedge product :
$ w' \otimes \bar{w} \mapsto  2 w' \wedge \bar{w} = w' \otimes \bar{w} - \bar{w}
\otimes w'$.

     Consider the  alternating bilinear form
$$A\in \Lambda^2(\Gamma
\otimes \C)^{\vee} \otimes (\Lambda \otimes \C ) = \Lambda^2(V \oplus
\overline{V})^{\vee} \otimes (U \oplus
\overline{U}) ,$$  satisfying the first bilinear relation and let us write
$$A = B + \overline{B},$$ where
$B \in \Lambda^2(\Gamma
\otimes \C )^{\vee}  \otimes U,$
and
$\overline{B} \in \Lambda^2(\Gamma
\otimes \C )^{\vee}  \otimes \overline{U}.$

By the first bilinear relation $B = B' + B''$, with $B' \in \Lambda^2(V)^{\vee}
\otimes U$, $B'' \in (V^{\vee} \otimes \overline{V})^{\vee} \otimes U.$

Concretely, $ A = B' + B'' + \overline { B'} +\overline { B'' }$, where
$B'$ is an alternating complex bilinear form.
The fact that $B''$ is
alternating reads out as:
$$B'' (v', \bar{v}) = - B'' (\bar{v}, v')  \   \forall v, v' \in V$$
whereas conjugation of tensors reads out as:
$$\bar{B} (\bar{x}, \bar{y}) = \overline {B (x, y)}\  \forall x,y \ \Rightarrow
\bar{B''} (\bar{v}, v') = \overline {B'' (v, \bar{v'})}.$$

\section{ Appell Humbert families}

We shall recall in this section the definition of the family given by the pairs of subspaces
satisfying the Riemann bilinear relations and some related results obtained in
\cite{cat04}, \cite{cf}.

\begin{df}
Given $A$ as above, we define $\mathcal
T \mathcal B_A$ as the subset of the product of
Grassmann Manifolds $Gr(m,2m)\times Gr(d,2d)$ defined by

$${\mathcal T}{\mathcal B}_{A}= \{(V,U) \in Gr(m,2m) \times Gr(d,2d)
\ | \ V \cap
\overline{V} = (0),$$
$$U \cap \overline{U} = (0), \ | \ {\rm the \ component \ of}  \ A \ in \
\Lambda^2(\overline{V})^{\vee} \otimes U \ is \ =0 \}.$$
This complex space is called the Appell Humbert space of Torus Bundles.

Since this space is not connected, we restrict ourselves to its
intersection with $\sT_m \times \sT_d$, i.e., we take two 
fixed orientations
of $\La$, resp. $\Ga$, and consider pairs of  complex
structures which have the same orientation as  the fixed ones.
\end{df}

One sees immediately that ${\mathcal T}{\mathcal B}_{A}$ is a complex
analytic variety of codimension at most $dm(m-1)/2$.

Note however that, for $ d \geq 3 , m > > 0$ we get a negative expected
dimension.  The structure of these complex spaces should be investigated in
general, for our present purposes we recall from \cite{cf} the
following

\begin{lem}
(\cite{cf})
If $ d=1$ then the open set ${\mathcal T}{\mathcal B}_{A}
\cap (\sT_m \times \sT_d)$ is connected.
\end{lem}

\begin{df}\label{ah}

The standard (Appell-Humbert) family of torus bundles parametrized by
${\mathcal T}{\mathcal B}_{A}$ is the family of principal holomorphic torus
bundles $X_{V,U}$ over $Y:=Y_V$ and with fibre $T:=T_U$
determined by the cocycle in $H^1(Y, {\mathcal H}(T)_{Y})$ obtained by taking
$f_{\gamma}(v)$ which is the class modulo $\Lambda$ of
$$F_{\gamma}(v):= B'(v, p_V(\ga)) + 2 B''(v, \overline{p_V (\ga)})
     + B''(p_V (\ga), \overline{p_V (\ga)}),\forall v \in V.$$

In other words,  $X_{V,U}$ is the quotient of $  T_U  \times V$
by the action of $\Ga$ such that
$$\ga  ( [u], v)  = (  [u +  B'(v, p_V(\ga)) + 2 B''(v, \overline{p_V (\ga)})
     + B''(p_V (\ga), \overline{p_V (\ga)})], v + p_V (\ga)).$$

\end{df}

\begin{oss}
As observed in \cite{cf} the above formula was a correction of the formula given in Definition 6.4 of
\cite{cat04}, where an identification of $\Ga \otimes \R$ with $V$ was
used, and thus $ A(z, \ga)$ was identified with $ B (z, \ga)$.
In the latter formula one had thus

$ -B'(v, \ga) -  B''(v, \overline{p_V (\ga)})$ instead of

     $B'(v, \ga) + 2 B''(v, \overline{p_V (\ga)})
     + B''(p_V (\ga), \overline{p_V (\ga)})$.

\end{oss}

We also recall from \cite{cat04} the definition of the complete
Appell-Humbert space.

\begin{df}

Given $A$ as above we define

$${\mathcal T}'{\mathcal B}_{A} = \{(V,U, \phi) \ | \ (V,U)
\in {\mathcal T}{\mathcal B}_{A}, \
\phi \in H^1(Y_V, {\mathcal H}(U)_{Y_V}) \cong
\overline{V}^{\vee} \otimes U \}.$$

The complete Appell-Humbert family of torus bundles parametrized by
${\mathcal T}'{\mathcal B}_{A}$ is the family of principal holomorphic torus
     bundles $X_{V,U, \phi}$ on $Y := Y _V$ and with fibre
     $T : = T_U$ determined by the cocycle in $H^1(Y, {\mathcal H}(T)_{Y})$
obtained by taking the sum of $f_{\gamma}(z)$ with the cocycle
$\phi \in H^1(Y_V, {\mathcal H}(U)_{Y_V}) \cong H^1(Y,{\mathcal
O}^d_{Y})$.

\end{df}

Finally we have the following

\begin{teo}{\bf \cite{cat04}}

Any principal holomorphic torus bundle with extension class isomorphic to
$\epsilon \in H^2(\Gamma, \Lambda)$ occurs in the complete Appell-Humbert
family ${\mathcal T}'{\mathcal B}_{A}$.

\end{teo}

We recall from \cite{cf} also the following

\begin{prop}\label{smooth}
\cite{cf}
Let $A: \Gamma \times \Gamma \rightarrow \Lambda$ be non zero.

     If $m=2, d=1$, i.e., $\Gamma \cong {\Z}^4$, $\Lambda \cong {\Z}^2$, both
Appell - Humbert spaces
     ${\mathcal T}{\mathcal B}_{A}$ and
${\mathcal T}'{\mathcal B}_{A}$ are smooth.

If $d =1$ and $m \geq 3$, ${\mathcal T}{\mathcal B}_{A}$
     is singular at the points where
$B'' =0$.

\end{prop}

\section{Symmetries of principal holomorphic torus bundles over tori}

Let us now assume that $f: X \rightarrow Y$ is a real principal holomorphic
torus bundle over a torus, and let $\sigma: X \rightarrow X$ be an
antiholomorphic involution. We assume as above that $Y = Y_V$, $V \cong
{\C}^m$, while the fibre of $f$ is $T = T_U$, with $U \cong {\C}^d$. 
In \ref{unicover} we have seen that the universal covering of $X$ is 
isomorphic to $V \oplus U \cong {\C}^{m+d}$ and we can find a lifting $\tilde{\sigma}$ of
$\sigma$ to the universal covering $U \oplus V\cong {\C}^{m+d}$. 

\begin{prop}
\label{restru}
Assume that the alternating form $A: \Gamma \times \Gamma \rightarrow \Lambda$ 
is nondegenerate, or equivalently that $\Lambda = Z(\Pi)$.
Let $\tilde{\sigma}$ be a lifting of $\sigma$ to $U \oplus V \cong {\C}^{m+d}$:
 then $
\tilde{\sigma}$ is an affine transformation. 
\end{prop}
{\bf Proof.}
Since $\tilde{\sigma}$ is a lifting of $\sigma$ to the universal covering, $\tilde{\sigma}$ 
acts by conjugation on the fundamental group $\Pi$ of $X$, therefore it also acts on the centre
$\Lambda$ of $\Pi$, because it is characteristic. 
 
Hence for every $\lambda \in \Lambda$, there exists a $\lambda' \in \Lambda$ such that 
$$\tilde{\sigma}(u +p_U(\lambda),v) = \tilde{\sigma}(u,v) + p_U(\lambda'),$$
where $p_U: U \oplus \bar{U} \rightarrow U$ denotes as usual the projection on
the first factor.
Assume that $\tilde{\sigma} : U \oplus V \rightarrow U \oplus V$
is given by
 $\tilde{\sigma}(u,v) = (\sigma_1(u,v), \sigma_2(u,v))$. 
Then we must have 
$$(\sigma_1(u+ p_U(\lambda),v) , \sigma_2(u+ p_U(\lambda),v)) = 
(\sigma_1(u,v) + p_U(\lambda'), \sigma_2(u,v)),$$
therefore $\sigma_2$ is constant as a function of $u$, 
and we may write $\sigma_2(u,v) = \sigma_2(v)$. 

Looking at the first component we obtain $\sigma_1(u + p_U(\lambda),v) = 
\sigma_1(u,v) + p_U(\lambda')$, so $\sigma_1$ is affine antiholomorphic in $u$ 
and we may write 
\begin{equation}
\label{sigma1}
\sigma_1(u,v) = A_1(v) \bar{u} + c_1(v)
\end{equation}
where $A_1(v)$ is a  linear map depending antiholomorphically
on $v$ (we may think of it as a $(d \times d)$ matrix once we fix a basis for $U$).

Recall that, for every $\gamma \in \Gamma$, and for any lift to $U \oplus V$
of the action of $\ga$   on $ T \times V$, there exists $\lambda'
\in \Lambda$ such that 
$${\gamma}(u,v) = (u + F_{\gamma}(v) + p_U(\lambda'), v +
p_V(\gamma)), \ \forall (u,v) \in U \oplus V,$$ where $F_{\gamma}(v)$ is as in definition
\eqref{ah} and $p_V: V \oplus \bar{V} \rightarrow V$ is the
first projection.
For all $\gamma \in \Gamma$ there must therefore 
exist $\gamma' \in \Gamma$, and a $\lambda''
\in \Lambda$ such that $\tilde{\sigma} \circ {\gamma} = \lambda'' {\gamma'} \circ
\tilde{\sigma}$.  Hence we have 
$$\tilde{\sigma} \circ {\gamma}(u,v) = 
\tilde{\sigma}(u +F_{\gamma}(v) + p_U(\lambda'),v + p_V(\gamma) )= $$
$$=(\sigma_1(u
+F_{\gamma}(v) + p_U(\lambda') ,v + p_V(\gamma)), \sigma_2(v + p_V(\gamma)))= $$
$$
=\lambda''{\gamma'}(\tilde{\sigma}(u,v)) = (\sigma_1(u,v)
+F_{\gamma'}(\sigma_2(v)) + p_U(\lambda'''), \sigma_2(v) + p_V(\gamma')),$$
where $\lambda''' \in \La$.
Therefore we obtain that $\sigma_2$ is affine antiholomorphic, 
$$\sigma_2(v) = A_2 \bar{v} + d_2,$$
where $A_2$ is a linear map (a $(m \times m)$ matrix if we fix a basis of $V$)
and we have 
$$A_2 \bar{v} + A_2 \overline{p_V(\gamma)} = A_2 \bar{v} + p_V(\gamma'),$$ 
so  that
$$A_2 \overline{p_V(\gamma)} = p_V(\gamma').$$
Looking at the first component we have 
$$\sigma_1(u + F_{\gamma}(v)+ p_U(\la'), v + p_V(\gamma)) = \sigma_1(u,v) +
 F_{\gamma'}(A_2 \bar{v} +d_2) + p_U(\lambda''').$$
Now using (\ref{sigma1}) we have 
\begin{equation}
\label{c1}
A_1(v + p_V(\gamma))( \bar{u} + \overline{F_{\gamma}(v)} +
\overline{p_U(\la')}) + c_1(v + p_V(\gamma)) = 
A_1(v) \bar{u} + c_1(v) + F_{\gamma'}(A_2 \bar{v}) + F_{\gamma'}(d_2) + p_U(\lambda''').
\end{equation}

So by derivation with respect to the variables $\bar{u_i}'s$ we obtain 
$$A_1(v + p_V(\gamma)) = A_1(v),$$
for all $v \in V$, for all $\gamma \in \Gamma$, therefore $A_1$ is constant in $v$ 
and we can write $A_1(v) = A_1$.
Now (\ref{c1}) becomes
\begin{equation}
\label{c1dinuovo}
c_1(v + p_V(\gamma)) -c_1(v) = F_{A_2 \overline{p_V(\gamma)}}(A_2 \bar{v}) - 
A_1 \overline{F_{\gamma}(v)} + F_{{A_2\overline{p_V(\gamma)}}}(d_2) +
 p_U(\lambda''')- A_1 \overline{p_U(\la')}
\end{equation}
and derivation with respect to the variables $\bar{v}_j's$ yields
the vanishing of the derivatives of  
$$c_1(v +p_V(\gamma)) -c_1(v) , \ \forall v \in V, \ \forall \gamma \in \Gamma,$$
so  these derivatives are constant and 
$$c_1(v) = q(v,v) + l(v) + d_1,$$ 
where $q(v,v)$ is quadratic in $v$, $l(v)$ is linear in $v$, and $d_1$ is a constant.

Now (\ref{c1dinuovo}) gives 
$$
\begin{gathered}
q(v + p_V(\gamma), v + p_V(\gamma)) + l(v)  + l(p_V(\gamma)) + d_1 - q(v,v) -l(v)
-d_1 = \\
=F_{{A_2 \overline{p_V(\gamma)}}}(A_2 \bar{v}) - A_1 \overline{{F_{\gamma}(v)}}
+ F_{{A_2 \overline{p_V(\gamma)}}}(d_2) + p_U(\lambda''')- A_1 \overline{p_U(\la')},
\end{gathered}
$$
so
$$q(p_V(\gamma), p_V(\gamma)) + 2q(v, p_V(\gamma)) + l(p_V(\gamma)) = $$
$$ =F_{{A_2 \overline{p_V(\gamma)}}}(A_2 \bar{v}) - A_1 \overline{{F_{\gamma}(v)}}
 + F_{{A_2 \overline{p_V(\gamma)}}}(d_2) + p_U(\lambda''')- A_1 \overline{p_U(\la')}.$$
Now, by looking in the above expression at the variable $\gamma$ we immediately get
$q(p_V(\gamma), p_V(\gamma)) = 0$, since it is the only quadratic term 
(substitute $\ga$ with $m \ga$ 
 and look at the asymptotic growth).  Since $q(p_V(\gamma), p_V(\gamma)) =
0$ for all $\gamma \in \Gamma$ and $p_V(\Gamma)$ is a lattice in $V$, 
we must have $q =0$ and
$c_1(v) = l(v) + d_1 = L \bar{v} + d_1$, where $L$ is a $(d \times m)$ matrix. 
So we finally get 
$$\tilde{\sigma}(u,v) = (A_1 \bar{u} + L\bar{v} + d_1, A_2 \bar{v} + d_2)$$
and the proposition is proven.
\hfill{Q.E.D.}

\begin{REM}
With the above notation we have 
$$\tilde{\sigma}(u,v)=  
\left( \begin{array}{cc} 
A_1 & L \\ 
0 & A_2
\end{array} \right)
\left( \begin{array}{c} 
\bar{u}\\ 
\bar{v}
\end{array} \right)
+ \left( \begin{array}{c} 
d_1\\ 
d_2
\end{array} \right)
$$
and the following properties hold:
\begin{enumerate}
\item $A_1 \overline{p_U(\Lambda)} = p_U(\Lambda)$, $A_2
  \overline{p_V(\Gamma)} = p_V(\Gamma)$. 
\item $\forall \ga \in \Ga$ we have
$$A_1 \ \overline{B''(p_V(\ga), \overline{p_V(\ga)})} = B''(A_2 \overline{p_V(\ga)}, 
\overline{A_2} p_V(\ga)).$$
\item $\forall v \in V$, $\forall \ga \in \Ga$ we have 
$$A_1 \overline{B'(v, p_V(\ga))} + 2 A_1 \ \overline{B''(v,
  \overline{p_V(\ga)})}  = B'(A_2\overline{v}, A_2 \overline{p_V(\ga)}) +
 2B''(A_2\overline{v}, \overline{A_2} p_V(\ga)) .$$
\item $\forall \ga \in \Ga$ 
$$L \overline{p_V(\ga)} - B'(d_2, A_2 \overline{p_V(\ga)}) + 2 B''(d_2,
\overline{A_2} p_V(\ga)) + B''(A_2 \overline{p_V(\ga)}, \overline{A_2}
p_V(\ga)) \in p_U(\La).$$
\item $A_1 \bar{A_1} = I$,  $A_2 \bar{A_2} = I$.
\item There exists a $\ga \in \Ga$ such that $\forall v \in V$ 
$$A_1 \overline{L}(v) + L \overline{A_2} v = F_{\ga}(v).$$
\item $\sigma_2^2 = Id$  mod $p_V(\Gamma)$, i.e. $\sigma_2$ induces an
  antiholomorphic involution on $Y$. 
\item $A_1 \overline{d_1} + L \overline{d_2} + d_1 \in p_U(\La)$. 
 
\end{enumerate}

\end{REM}  
{\bf Proof.}
Conditions (1), (2), (3) and (4) easily follow as in the proof of
\eqref{restru} by imposing that for all $\la
\in \La$ there exists a $\la' \in \La$ such that $\tilde{\sigma} \circ \la=
\la' \circ \tilde{\sigma}$, and for all $\gamma \in \Ga$ there must exist
$\gamma' \in \Ga$, $\lambda'' \in \La$ such that $\tilde{\sigma} \circ \ga=
\la''\ga'  \circ \tilde{\sigma}$.
Here we also used the expression of $F_{\ga}(v)$ given in definition
\eqref{ah}:
$$F_{\gamma}(v):= B'(v, p_V(\ga)) + 2 B''(v, \overline{p_V (\ga)})
     + B''(p_V (\ga), \overline{p_V (\ga)}),\forall v \in V.$$

Conditions (5), (6), (7) and (8) immediately follow by imposing $\tilde{\sigma}^2 \in
\Pi$.
\qed

\begin{teo}
Let $f:X \rightarrow Y$ be a principal holomorphic torus bundle over a torus
such that $A$ is non degenerate and assume that $\sigma: X \rightarrow X$ is
an antiholomorphic involution  on $X$. The differentiable type of the pair $(X, \sigma)$ 
is completely determined by the orbifold fundamental group exact sequence. 
More precisely, the affine embedding of $\Pi_{\sigma}$ is uniquely determined
up to conjugation.
\end{teo}
{\bf Proof.}
Let 
\begin{equation}
\label{orfun}
1 \rightarrow \Pi \rightarrow \hat{\Pi} := \Pi_{\sigma} \rightarrow {\Z}/2 \rightarrow 1
\end{equation}
be the orbifold fundamental group exact sequence of the pair $(X, \sigma)$. 
Every lifting $\tilde{\sigma}$ of $\sigma$ to $\hat{\Pi}$ acts by conjugation
 on $\Pi$, so it acts by
conjugation on the centre $\Lambda$ of $\Pi$ and thus it acts on the quotient $\Gamma =
\Pi/\Lambda$. Therefore we have determined an extension 
$$1 \rightarrow \Gamma \rightarrow \hat{\Gamma} \rightarrow {\Z}/2 \rightarrow 1$$
which is the orbifold fundamental group exact sequence of the real torus $(Y, \sigma_2)$,
 where $\sigma_2$ denotes as above the second component of $\tilde{\sigma}$.   
Since for a real torus the orbifold fundamental group exact sequence determines 
the differentiable type (cf. \cite{cat02}), we have shown that we can fix the 
differentiable type of the pair
$(Y ,\sigma_2)$.  

By proposition \ref{restru} we know that any lifting $\tilde{\sigma}$ of
$\sigma$ to the universal covering $(\Lambda \otimes {\R}) \oplus
 (\Gamma \otimes {\R})$ is of the form

$$\tilde{\sigma}(y,x)=  
\left( \begin{array}{cc} 
A_1 & L \\ 
0 & A_2
\end{array} \right)
\left( \begin{array}{c} 
y\\ 
x
\end{array} \right)
+ \left( \begin{array}{c} 
d_1\\ 
d_2
\end{array} \right)
$$
and we know the $(2m \times 2m)$-matrix $A_2$ and the translation vector $d_2$.
We also have $A_1^2 = I$, $A_2^2 = I$, because $\tilde{\sigma}^2 \in \Pi$. 
Furthermore we know $\tilde{\sigma}^2 \in \Pi$, since we know the extension (\ref{orfun}),
therefore we know the vector   
$$ 
\left( \begin{array}{cc} 
A_1 & L \\ 
0 & A_2
\end{array} \right)
\left( \begin{array}{c} 
d_1\\ 
d_2
\end{array} \right)
+ \left( \begin{array}{c} 
d_1\\ 
d_2
\end{array} \right)=
\left( \begin{array}{c}  
L d_2 + A_1 d_1 + d_1\\
A_2 d_2 + d_2
\end{array} \right)
$$
which is the translation part of $\tilde{\sigma}^2$.

For every $\lambda \in \Lambda$ we know $\tilde{\sigma} \lambda 
{\tilde{\sigma}}^{-1} \in \Lambda$ and $\tilde{\sigma} 
\lambda {\tilde{\sigma}}^{-1}(w) = w +  M
\lambda$, $\forall w \in (\Gamma \otimes {\R}) \oplus (\Lambda \otimes {\R})$, where 
$$M =  
\left( \begin{array}{cc} 
A_1 & L \\ 
0 & A_2

\end{array} \right)
$$
Thus we know $A_1(\lambda) \in \Lambda$, $\forall \lambda \in \Lambda$ and since
 $\Lambda$ generates $\Lambda \otimes {\R}$, we know $A_1$.
So we know $A_1$, $A_2$, $d_2$, and $L d_2 + A_1 d_1 + d_1$.

For all $\hat{\gamma} \in \Pi$ lifting  a given $\gamma \in \Gamma$, $\gamma \neq 0$, 
 we know $\tilde{\sigma} \hat{\gamma} \tilde{\sigma}^{-1}$. Let us set $\tilde{\sigma}(w) = 
Mw +
b$, $\hat{\gamma}(w) = Dw + h$, where 
$$D \left( \begin{array}{c} 
y\\ 
x
\end{array} \right)=  
\left( \begin{array}{cc} 
I_{2d} & \phi_{\gamma}\\ 
0 & I_{2m}
\end{array} \right)
\left( \begin{array}{c} 
y\\ 
x
\end{array} \right)
$$
where $\phi_{\ga}(x)$ can be chosen to be equal to $A(x, \ga)$ as it is proven
in \eqref{difaction1}, 
and $h = \left( \begin{array}{c} 
l\\ 
\gamma
\end{array} \right)$, where $l \in \Lambda$.

We have $\tilde{\sigma} \hat{\gamma} \tilde{\sigma}^{-1}(w) = \tilde{\sigma}
 \hat{\gamma} (M^{-1}w - M^{-1}b) = \tilde{\sigma}(DM^{-1}w -DM^{-1}b +h) = MDM^{-1} w -
MDM^{-1} b + Mh + b$. Thus we know $MDM^{-1}$ and $-MDM^{-1}b + Mh +b$. One easily
computes 
$$-MDM^{-1}b + Mh +b = \left( \begin{array}{c} 
-A_1 \phi_{\gamma'}A_2^{-1}d_2 + L \gamma + A_1 l\\
A_2 \gamma
\end{array} \right)
$$
and since we know both $-A_1 \phi_{\gamma'}A_2^{-1}d_2$ and  $A_1 l$, we also
know $L \gamma$ for all $\gamma \in \Gamma$. Now $\Gamma$ generates $\Gamma
\otimes {\R}$, so we know $L: \Gamma \otimes {\R} \rightarrow \Lambda \otimes
{\R}$.

We have already seen that we know $L d_2 + A_1 d_1 + d_1$, so we also know
$d_1' = A_1 d_1 + d_1$ and we know $A_1$. 
Now, $A_1^2 = I$, so we can decompose $\Lambda \otimes {\R} = Z \oplus S
\oplus W^+ \oplus W^-$, where $W^{\pm}$ are the $\pm 1$-eigenspaces of $A_1$
and where 
$$A_1(z,s,w^+, w^-) = (s,z,w^+, -w^-).$$
Observe that since $A_1$ is antiholomorphic, $W^+$ and $W^-$ have
 the same dimension (also $dim (Z) = dim (S)$). Thus we can write $d_1 = (z,s,d_1^+, d_1^-)$.
$A_1(z,s,d_1^+, d_1^-) + (z,s,d_1^+, d_1^-)= (s+z,s+z, 2d_1^+, 0)$, so we know
$d_1^+$ and $z+s$. We can change the origin by translating with $(y,x) \in
(\Lambda \otimes {\R}) \oplus (\Gamma \otimes {\R}) \mapsto (y+t, x) 
\in (\Lambda \otimes {\R}) \oplus (\Gamma \otimes {\R})$ and we may assume
 that $d_1 = (z , s ,
d_1^+, d_1^-) + A_1(t) -t = (z, s, d_1^+, d_1^-) + (t_2-t_1, t_1-t_2, 0, -2t^-)$ 
($t = (t_1, t_2, t^+,
t^-)$). So we can choose $t^- = d_1^-/2$,  $t_2 - t_1 = -z$, therefore the first and the last
components of $d_1$ can be chosen equal to zero. Then the second component is 
$s +  t_1-t_2 = s+
z $ and therefore we know it, finally the third component is $d_1^+$ and we already know it.
 
\qed

\begin{REM}
\label{top}
Let us fix the orbifold fundamental group exact sequence of a real principal holomorphic torus
 bundle over a torus such that the alternating bilinear form $A: \Gamma \times \Gamma
\rightarrow \Lambda$ is non degenerate.  Then  the topological and differentiable
action of $\Pi_{\sigma}$ is fixed, and the action of any element $\tilde{\sigma}$ of the orbifold
fundamental group
$\hat{\Pi}$ on the universal covering $(\La \otimes {\R}) \oplus  (\Ga \otimes {\R})$ 
is  given by the
affine transformation 
$$\tilde{\sigma}
\left( \begin{array}{c} 
y\\ 
x
\end{array} \right)
=
\left( \begin{array}{cc} 
A_1 & L \\ 
0 & A_2
\end{array} \right)
\left( \begin{array}{c} 
y\\ 
x
\end{array} \right)
+ 
\left( \begin{array}{c} 
d_1\\ 
d_2
\end{array} \right)
$$
If we fix such an element   $\tilde{\sigma}$, the following holds:

\begin{enumerate}
\item $A_1 : \La \rightarrow \La$, $A_1^2 = I$. 
\item $A_2: \Ga \rightarrow \Ga$, $A_2^2 = I$.
\item $\forall \gamma, \ \forall x \in \Gamma \otimes \R$, $A_1(A(x, \gamma)) = 
A(A_2(x), A_2(\gamma))$.
\item The ${\R}$ - linear map $L': \Gamma \otimes {\R} \rightarrow \Lambda \otimes 
{\R}$, $L'(x) : = Lx - A(d_2, A_2 (x))$ satisfies $L'(\Gamma) \subset \Lambda$.
\item $A_2(d_2) + d_2 \in \Gamma$.
\item $L(d_2) + A_1(d_1) + d_1 \in \Lambda$.
\item $ \exists \gamma \in \Gamma$ such that $L(A_2(x)) + A_1(L(x)) = -A(x,
  \gamma)$ $\forall x \in \Gamma \otimes \R$. 
\end{enumerate}
\end{REM}
{\bf Proof.}
For all $\lambda \in \Lambda$ we have $\tilde{\sigma} \lambda \tilde{\sigma}^{-1} 
\in \Lambda$ and this immediately implies that $A_1(\Lambda) \subset \Lambda$. $A_1^2 = I$,
since $\tilde{\sigma}^2 \in \Pi$ and for all $g \in \Pi$ the action of $g$ on $(\La \otimes {\R})
\oplus (\Ga \otimes {\R})$ is given by
$$g
\left( \begin{array}{c} 
y\\ 
x
\end{array} \right)
=
\left( \begin{array}{cc} 

I & \phi_{\ga}\\ 
0 & I
\end{array} \right)
\left( \begin{array}{c} 
y\\ 
x
\end{array} \right)
+ 
\left( \begin{array}{c} 
l\\ 
\ga
\end{array} \right)
$$
where $\gamma$ is the image of $g$ in $\Gamma$, $l \in \Lambda$,  and
$\phi_{\ga}(x) = A(x, \ga)$ as in \eqref{difaction1}. 

This also implies that $A_2^2 = Id$. 

For every $g \in \Pi$, $g \not\in \Lambda$, we know $\tilde{\sigma} g \tilde{\sigma}^{-1} 
= h \in \Pi - \Lambda$. If we set 
$$g
\left( \begin{array}{c} 
y\\ 
x
\end{array} \right)
=
\left( \begin{array}{cc} 
I & \phi_{\ga} \\ 
0 & I
\end{array} \right)
\left( \begin{array}{c} 
y\\ 
x
\end{array} \right)
+ 
\left( \begin{array}{c} 
l\\ 
\ga
\end{array} \right)
$$
where as above $\gamma$ is the image of $g$ in $\Gamma$ and $l \in \Lambda$, and 
$$h
\left( \begin{array}{c} 
y\\ 
x
\end{array} \right)
=
\left( \begin{array}{cc} 
I & \phi_{\delta} \\ 
0 & I
\end{array} \right)
\left( \begin{array}{c} 
y\\ 
x
\end{array} \right)
+ 
\left( \begin{array}{c} 
\la\\ 
\delta
\end{array} \right)
$$
where $\delta$ is the image of $h$ in $\Gamma$ and $\lambda \in \Lambda$, we have 
$$\tilde{\sigma} g \tilde{\sigma}^{-1}
\left( \begin{array}{c} 
y\\ 
x
\end{array} \right)
=
\left( \begin{array}{cc} 
I & A_1 \phi_{\ga} A_2^{-1} \\ 
0 & I
\end{array} \right)
\left( \begin{array}{c} 
y\\ 
x
\end{array} \right)
+ 
\left( \begin{array}{c} 
-A_1 \phi_{\gamma}A_2^{-1} d_2 + L \gamma + A_1 l\\
A_2 \gamma
\end{array} \right)$$
$$
= h
\left( \begin{array}{c} 
y\\ 
x
\end{array} \right)
=
\left( \begin{array}{cc} 
I & \phi_{\delta}\\ 
0 & I
\end{array} \right)
\left( \begin{array}{c} 
y\\ 
x
\end{array} \right)
+ 
\left( \begin{array}{c} 
\la \\ 
\delta
\end{array} \right)
$$

This yields $A_1 \phi_{\gamma}A_2^{-1} = \phi_{\delta}$, $A_2 \gamma = \delta$,
 $-A_1 \phi_{\gamma} A_2^{-1} d_2 + L \gamma + A_1 l = \lambda \in \Lambda$. So
$A_2(\Gamma) = \Gamma$ and for all $\gamma \in \Gamma$ we have $L' \gamma = 
L \gamma -
\phi_{A_2 \gamma} (d_2) = L(\gamma) - A(d_2, A_2(\gamma)) \in \Lambda$. 

The condition $A_1 \phi_{\gamma} A_2^{-1} = \phi_{A_2(\gamma)}$, $\forall
\gamma \in \Gamma$ can be written as 
$$A_1(A(x, \gamma)) = A(A_2(x),A_2(\gamma)), \ \forall x \in \Ga \otimes \R, 
\ \forall \gamma \in \Gamma.$$

Finally $\tilde{\sigma}^2 \in \Pi$ implies that $A_2(d_2) + d_2 \in \Gamma$,
$L(d_2) + A_1(d_1) + d_1 \in \Lambda$, furthermore there must exist a $\gamma
\in \Gamma$ such that $L A_2 + A_1 L = \phi_{\gamma}$.
This implies that, for every $\delta \in \Gamma$, $L(A_2(\delta)) + A_1(L(\delta)) = 
-A(\delta, -\gamma) \in \Lambda$. 
\qed

We have three ${\R}$ - linear maps $A_1: \La \rightarrow \La$, $A_2: \Ga
\rightarrow \Ga$, $L : \Gamma \otimes \R  \rightarrow \Lambda \otimes \R$ with
the above properties. 

Let us now fix a complex structure on the bundle $f: X \rightarrow Y$, i.e. we
fix a point $(V, U)$ in ${\mathcal T}{\mathcal B}_A = \{(V, U) \in Gr(m,2m) \times Gr(d,2d) 
\ | \ V \cap \overline{V} = (0),
\ U \cap \overline{U} = (0), \ | \ the \
 component \ of \ A \ in \ \Lambda^2(\overline{V})^{\vee} \otimes U \ is \ =0 \}$. 

We may now see $A_1$ as a real element in $(\La \otimes {\C})^{\vee} \otimes (\La \otimes {\C})
 = (U \oplus \overline{U})^{\vee} \otimes (U \oplus \overline{U})$ and we want to impose that
$A_1$ is antiholomorphic with respect to the chosen complex structure, so we want that the
component of $A_1$ in $U^{\vee} \otimes U$ is zero. Analogously we see $A_2$ as a real
 element in
$(\Ga \otimes {\C})^{\vee} \otimes (\Ga \otimes {\C}) = (V \oplus \overline{V})^{\vee} \otimes (V
\oplus \overline{V})$ and we want that the component of $A_2$ in $V^{\vee} \otimes V$ is zero.
Finally we also see $L$ as a real element in $(\Gamma \otimes {\C})^{\vee} \otimes (\Lambda
\otimes {\C}) = (V \oplus \overline{V})^{\vee} \otimes (U \oplus \overline{U})$ 
and we want that
its component in $V^{\vee} \otimes U$ is zero.

Observe that since any other lifting $\tilde{\sigma}'$ of $\sigma$ to $\hat{\Pi}$ is
 obtained from $\tilde{\sigma}$ by composition with an element in $\Pi$, which acts
holomorphically with respect to the chosen complex structure, we may give the following
 definition.

\begin{df}

 We define spaces

$${\mathcal T}{\mathcal B}_{A}^{\R}(A_1,A_2,L)= \{(V,U) \in Gr(m,2m) \times Gr(d,2d) 
\ | \ V \cap \overline{V} = (0),$$
$$U \cap \overline{U} = (0), \ | \ the \
 component \ of \ A \ in \ \Lambda^2(\overline{V})^{\vee} \otimes U \ is \ =0, \  
the \ component$$
 $$of \ A_1 \ in \ U^{\vee} \otimes U \ is \ =0, \ the \ component \ of \ A_2 \ in 
\ V^{\vee} \otimes V \ is \ =0,$$
$$the \ component \ of \ L \ in \ V^{\vee} \otimes U \ is \ =0 \},$$  
and 
$${\mathcal T}'{\mathcal B}_{A}^{\R}(A_1,A_2,L) = \{(V,U, \phi) \ | \ (V,U)
\in {\mathcal T}{\mathcal B}_{A}^{\R}(A_1,A_2,L) \
\phi \in H^1(Y_V, {\mathcal H}(U)_{Y_V}) \cong
\overline{V}^{\vee} \otimes U \}.$$

\end{df}

\begin{oss}
Assume that we have fixed the orbifold fundamental
group exact sequence of the pair $(X, \sigma)$ and that the alternating
bilinear form  
$A: \Gamma \times \Gamma \rightarrow \Lambda$ is non degenerate:
then the above spaces parametrize families of real structures
on   principal holomorphic torus bundles, and in particular 
the latter family parametrizes all the possible real structures on a  
principal holomorphic torus bundle
over a torus with a  given real topological type.
\end{oss}

Let us now write 
$$\La \otimes \R = U^+ \oplus U^-,$$
$$\Ga \otimes \R = V^+ \oplus V^-,$$
according to the eigenspace decomposition  for $A_1$, respectively for $A_2$. 

If we set $A = A^+ + A^-$, where $A^+$ is the component of $A$ 
with values in $U^+$,
$A^-$ is the component of $A$  with values in $U^-$, we can write condition (3) of
(\ref{top}) as follows: $\forall x \in \Ga \otimes \R$, $\forall \ga \in \Ga$, 
$$A_1(A(x, \ga)) = A^+(x, \ga) - A^-(x, \ga) = A(A_2(x), A_2(\ga)) = A(x^+ -
x^-, \ga^+ - \ga^-),$$
where $x = x^+ + x^-$, $\ga = \ga^+ + \ga^-$ according to the above
decomposition.

Thus we have 
$$A^+( x^+ + x^-, \ga^+ + \ga^-) = A^+( x^+ - x^-, \ga^+ - \ga^-), $$
$$A^-( x^+ + x^-, \ga^+ + \ga^-) = -A^-( x^+ - x^-, \ga^+ - \ga^-), $$
which in turn is equivalent to 
$$A^+( x^+, \ga^-) + A^+( x^-, \ga^+) = 0, \  \ A^-( x^+, \ga^+) + A^-( x^-, \ga^-) = 0.$$

We conclude then that condition (3) is equivalent to
\begin{equation}
\label{3.1}
A^+|_{ V^+ \times V^-} \equiv 0, \ A^-|_{ V^+ \times V^+} \equiv 0,  \ A^-|_{ V^- \times V^-}
\equiv 0 , \ i.e., 
\end{equation}
\begin{equation}
\label{3.2}
 A (x_1, x_2) = A^+ (x_1^+, x_2^+) + A^+ (x_1^-, x_2^-) + A^- (x_1^-, x_2^+) + A^- (x_1^+,
x_2^-).
\end{equation}

Condition (7) of (\ref{top}) is equivalent to the existence of $\hat{\ga} \in
\Ga$ such that 

$$L^+(x^+ - x^-) + L^-(x^+ - x^-) + L^+(x^+ + x^-) - L^-( x^+ + x^-) = $$
$$-A^+
(x^+ + x^-, \hat{\ga}^+ + \hat{\ga}^-) -A^-
(x^+ + x^-, \hat{\ga}^+ + \hat{\ga}^-),$$
that is, to

$$2 L^+(x^+) = - A^+(x^+, \hat{\ga}^+) - A^+(x^-, \hat{\ga}^-),$$ 
$$2 L^-(x^-) =  A^-(x^+, \hat{\ga}^-) + A^-(x^-, \hat{\ga}^+),$$
in particular we have
$$ A^+(x^-, \hat{\ga}^-)\equiv 0 \equiv A^-(x^+, \hat{\ga}^-),$$
or equivalently,  
$$A(-, \hat{\ga}^-) \equiv 0.$$
Since we are assuming $A$ nondegenerate, we must have $\hat{\ga}^- = 0$, and
we have 

\begin{equation}
\label{7}
\begin{gathered}
L^+(x^+) = - \frac{A^+}{2}(x^+, \hat{\ga}^+)\\
L^-(x^-) = \frac{A^-}{2}(x^-, \hat{\ga}^+)
\end{gathered}
\end{equation}

 We look  now for respective complex structures $J_1$, $J_2$  on $\La \otimes
\R$, $\Ga \otimes \R$ making $L$, $A_1$ and $A_2$ antiholomorphic. 

The condition $J_2 A_2 = -A_2 J_2$ immediately implies
that $J_2$ exchanges the two Eigenspaces for $A_2$, therefore we can write
 $J_2(x^+, x^-) = (Cx^-, D x^+)$.  The condition $J_2^2 = -Id$ is
then equivalent to $D = -C^{-1}$, so, if we set
$B_2 = -C$, we have 
$$J_2(x^+, x^-) = (-B_2x^-, B_2^{-1} x^+).$$
Proceeding analogously for $J_1$ we obtain 
$$J_1(y^+, y^-) = (-B_1y^-, B_1^{-1} y^+).$$

Finally $L \circ J_2 = -J_1 \circ L$ is equivalent to
$$L(-B_2x^-, B_2^{-1} x^+)= (B_1L^-(x), - B_1^{-1}L^+(x)),$$
equivalently,
$$L^+(-B_2x^-, B_2^{-1} x^+)= B_1L^-(x^+, x^-),$$
$$L^-(-B_2x^-, B_2^{-1} x^+)= -B_1^{-1}L^+(x^+, x^-).$$

If we now write $L(x^+, x^-)$ as $L_+ x^+ + L_- x^-$ we obtain 

$$-L^+_+B_2(x^-) + L^+_- B_2^{-1}(x^+) = B_1L^-_+(x^+) + B_1L^-_-(x^-),$$
$$-L^-_+B_2(x^-) + L^-_- B_2^{-1}(x^+)= -B_1^{-1} L^+_+(x^+) - B_1^{-1}L^+_-(x^-).$$

We rewrite the first equation as

\begin{equation}
\label{Lanti}
\begin{gathered}
L^+_+B_2 + B_1 L^-_- \equiv 0,\\
L^+_- B_2^{-1} - B_1 L^-_+ \equiv 0.
\end{gathered}
\end{equation}

After rewriting the second equation as
$$L^-_+ B_2 - B_1^{-1} L^+_- \equiv 0,$$
$$L^-_- B_2^{-1} + B_1^{-1} L^+_+ \equiv 0,$$
we observe that these equations are clearly equivalent to (\ref{Lanti}).

Conditions (\ref{7}) become 
\begin{equation} 
\label{7bis}
\begin{gathered}
L^+_+(x^+) = - \frac{A^+}{2}(x^+, \hat{\ga}^+),\\
L^-_-(x^-) =  \frac{A^-}{2}(x^-, \hat{\ga}^+).
\end{gathered}
\end{equation}

We shall now write the Riemann bilinear relations. 

Recall that if $(V,U) \in
{\mathcal T}{\mathcal B}_A$, $V = \{x - iJ_2 x \ | \ x \in \Ga \otimes \R\}$,
$U =  \{z - iJ_1 z \ | \ z \in \La \otimes \R\}$ and we want
$$A(x-iJ_2x, y-iJ_2y) = A(x,y) - A(J_2x, J_2y) -i(A(x, J_2y) + A(J_2x,
y)) \in U, \ \forall x, y \in \Ga \otimes \R.$$

This means that

\begin{equation}
\label{RBR}
A(x, J_2y) + A(J_2x, y) = J_1A(x,y) - J_1A(J_2x, J_2y).
\end{equation}
We shall now split (\ref{RBR}) into $U^+$ and $U^-$ components, so that using
(\ref{3.1}), (\ref{3.2}) the
first component is
\begin{equation}
\begin{gathered}
A^+_{++}(x^+, -B_2 y^-) + A^+_{--}(x^-, B_2^{-1}y^+) + \\
A^+_{++}(-B_2x^-, y^+)
+ A^+_{--}(B_2^{-1}x^+, y^-) = -B_1A^-_{+-}(x^+,y^-) - B_1A^-_{-+}(x^-,y^+)\\
- B_1A^-_{+-}(B_2 x^-,B_2^{-1}y^+)
- B_1A^-_{-+}(B_2^{-1} x^+, B_2 y^-)
\end{gathered}
\end{equation}
while the second component is 

\begin{equation}
\begin{gathered}
A^-_{+-}(x^+, B_2^{-1} y^+) + A^-_{-+}(B_2^{-1}x^+, y^+) + \\
A^-_{+-}(-B_2x^-, y^-)
+ A^-_{-+}(x^-, -B_2y^-) = B_1^{-1}[A^+_{++}(x^+,y^+) + A^+_{--}(x^-,y^-)\\
- A^+_{++}(B_2 x^-,B_2y^-)
- A^+_{--}(B_2^{-1} x^+, B_2^{-1} y^+)].
\end{gathered}
\end{equation}
From these equations we derive the following equations looking at the four
possible bilinear types $(x^+,y^+)$, $(x^+,y^-)$, $(x^-,y^+)$, $(x^-,y^-)$.

$$-A^+_{++}(x^+, B_2 y^-) + A^+_{--}(B_2^{-1}x^+, y^-) =$$
$$= -B_1A^-_{+-}(x^+,y^-) - B_1A^-_{-+}(B_2^{-1} x^+, B_2 y^-)$$

$$A^+_{--}(x^-, B_2^{-1}y^+) - A^+_{++}(B_2x^-, y^+)=$$
$$ =- B_1A^-_{-+}(x^-,y^+) - B_1A^-_{+-}(B_2 x^-,B_2^{-1}y^+)$$

$$A^-_{+-}(x^+, B_2^{-1} y^+) + A^-_{-+}(B_2^{-1}x^+, y^+) =$$
 $$= B_1^{-1}[A^+_{++}(x^+,y^+) - A^+_{--}(B_2^{-1} x^+, B_2^{-1} y^+)]$$

$$-A^-_{+-}(B_2x^-, y^-) - A^-_{-+}(x^-, B_2y^-) =$$
$$ = B_1^{-1}[A^+_{--}(x^-,y^-)
- A^+_{++}(B_2 x^-,B_2y^-)].$$

These four equations can be rewritten as tensor equations as follows
(according to the standard notation for the transformation of bilinear forms)
\begin{equation}
\label{1}
-A^+_{++}B_2 + ^tB_2^{-1} A^+_{--} = B_1 (-A^-_{+-} - ^tB_2^{-1}A^-_{-+}B_2)
\end{equation}

\begin{equation}
\label{2}
A^+_{--}B_2^{-1} - ^tB_2 A^+_{++} = B_1 (-A^-_{-+} - ^tB_2A^-_{+-}B_2^{-1})
\end{equation}

\begin{equation}
\label{3}
A^-_{+-}B_2^{-1} + ^tB_2^{-1} A^-_{-+} = B_1^{-1} (A^+_{++} - ^tB_2^{-1}A^+_{--}B_2^{-1})
\end{equation}

\begin{equation}
\label{4}
-^t B_2A^-_{+-} -A^-_{-+}B_2 = B_1^{-1} (A^+_{--} - ^tB_2A^+_{++}B_2)
\end{equation}

We observe that the tensor equations above are all equal, in fact (\ref{2})
yields (\ref{1}) by composing with $B_2$ to the right and
with  $^tB_2^{-1}$ to the left. (\ref{3})
yields (\ref{1}) by composing with $B_2$ to the right. (\ref{4})
yields (\ref{1}) by composing with  $^tB_2^{-1}$ to the left. 

So we have shown that if we have a real structure we have only one equation
for the Riemann bilinear relation:
\begin{equation}
\label{RBR1}
A^+_{--} - ^tB_2 A^+_{++} B_2 = -B_1(A^-_{-+}B_2 - ^tB_2 \ ^tA^-_{-+}),
\end{equation}
where we have used $A^-_{+-} = - ^tA^-_{-+}$ since $A^-$ is alternating.

Now, in order to simplify the notation, we set $A^+_{--} =: A_-$, $A^+_{++} =: A_+$,
$D:= A^-_{-+}$, so (\ref{RBR1}) becomes
\begin{equation}
\label{RBR2}
A_- - ^tB_2 A_+ B_2 = -B_1(DB_2 - ^t(DB_2)), 
\end{equation}

We can now easily show that in the case $m = d =1$ the space $ {\mathcal T}'{\mathcal
  B}_{A}^{\R}(A_1,A_2,L) \cap ({\mathcal T}_1 \times {\mathcal T}_1)$ is
  connected. We observe, here and in the following,  that it suffices to show that $ {\mathcal
T}{\mathcal
  B}_{A}^{\R}(A_1,A_2,L) \cap ({\mathcal T}_1 \times {\mathcal T}_1)$ is
  connected. Note that the case $m = d =1$, since we assume $A$ to be nondegenerate,
 corresponds to the case of Kodaira surfaces.
The fact
  that the moduli space of real Kodaira surfaces of a given topological type
  is connected was already proved by the second author with different methods
  (cf. \cite{f}).

\begin{prop}
If $m = d =1$, the space $ {\mathcal T}{\mathcal
  B}_{A}^{\R}(A_1,A_2,L) \cap ({\mathcal T}_1 \times {\mathcal T}_1)$ is
  connected.

\end{prop}

\Proof
We observe first of all that equation (\ref{RBR2}) does not appear since for
$m =1$, $A_- = A_+ = 0$, and $B_1$, $B_2$ are scalars. 

So we only have to consider conditions (\ref{Lanti}). With an appropriate
choice of orientation on each vector space we may assume $B_1>0, B_2>0$, so if
$L_-^- \neq 0$, by (\ref{Lanti}) we obtain $B_1 = - \frac{B_2 L^+_+}{L^-_-}$
and 
$$L^+_- + \frac{B_2^2 L^+_+ L^-_+}{L^-_-} =0,$$
or equivalently 
$$L^+_- L^-_- + B_2^2 L^+_+ L^-_+ =0,$$
which is solvable if either $L^+_- L^-_- = L^+_+ L^-_+ =0,$ or both are non zero
and have opposite sign and in this case we only have one positive solution.

If $L^-_- =0$, then also $L^+_+=0$ and we have $B_1B_2 =
\frac{L^+_-}{L^-_+}$. Since we are assuming $B_1>0$, $B_2>0$, we must have
$\frac{L^+_-}{L^-_+}>0$ and the set is clearly connected. 
\qed

\begin{REM}
\label{L=0}
Observe that if $L =0$ and $d =1$, we set $B_2 =:B$, $B_1 =:b$, a scalar,
and we only have the equation

\begin{equation}
A_- - ^tB A_+ B = -b(DB - ^t(DB)), 
\end{equation}

If we also assume $m =2$, we can write 
$A_- = a_-  \left( \begin{array}{cc} 
0& 1\\ 
-1 & 0
\end{array} \right)$, 

$A_+ = a_+  \left( \begin{array}{cc} 
0& 1\\ 
-1 & 0
\end{array} \right)$, 
and we obtain only one scalar equation 

$$a_- - det(B)a_+ = -b(d_{11} b_{12} + d_{12} b_{22} - d_{21} b_{11} - d_{22} b_{21})$$
and we must combine it with the condition 
$$det(B) = b_{11}b_{22} - b_{12} b_{21} >0.$$
\end{REM}

\begin{teo}
For $d =1$, $m =2$, the space ${\mathcal X} = {\mathcal
T}{\mathcal
  B}_{A}^{\R}(A_1,A_2,L) \cap ({\mathcal T}_2 \times {\mathcal T}_1)$ is connected.
\end{teo}

\Proof
By remark (\ref{L=0}), if $L=0$ we have the two conditions:

\begin{equation}
\label{l=0}
a_- - det(B)a_+ = -b(d_{11} b_{12} + d_{12} b_{22} - d_{21} b_{11} - d_{22}
b_{21})
\end{equation}

\begin{equation}
\label{detB}
det(B) = b_{11}b_{22} - b_{12} b_{21} >0,
\end{equation}
and we assume $b>0$. 

Assume that the matrix $D = 0$ : then we have a product of the half-line $ \{ b > 0 \}$
with the quadric $ \{ det B = a \}$, where $a$ is a positive constant. In this case we
are done since  $ \{ det B = a \}$ is a central quadric in $\R^4$ with quadratic
part of signature $(+,+,-,-,)$, and is therefore connected.

Assume then that $ D \neq 0$, then the rigth hand side is equal to $ b b'$,
where we define $ b'$ as the linear form 
$b' : = d_{11} b_{12} + d_{12} b_{22} - d_{21} b_{11} - d_{22}
b_{21}$. 

We want to change coordinates in the $\R^4$ of the matrices $B$, completing
the linear form $b'$ to a basis $ (b',x,y,z)$ of the space of linear forms.
The projective quadric  $ \{ det B = 0 \}$ and the hyperplane $ \{ b' = 0 \}$
determine an affine quadric in $\R^3$ with the given signature, therefore there are
only two cases, according to the property whether the hyperplane is transversal
or tangent to the quadric.

In the first case, after dehomogenizing (i.e., setting $b' = 1$)  we get
a  one sheeted (hyperbolic) hyperboloid:

$$ \pm det B =  {b'}^2 - z^2  + xy ,$$

in the second case we get a hyperbolic paraboloid

$$ det B =  {b'} z + xy.$$

Let us observe that 
the space  ${\mathcal X}$   maps to the open set $\Omega$ in $\R^4$
defined by the inequality $ det B > 0$, and that
if we define $ \Omega ' : = \Omega \cap \{ b' \neq 0 \}$,
then  $ \Omega ' $ is homeomorphic to its inverse image in ${\mathcal X}$.

It is immediate to conclude, in view of the above normal forms,
that the open sets  $ \Omega ^+ : = \Omega \cap \{ b' > 0 \}$,
resp. $ \Omega ^- : = \Omega \cap \{ b' < 0 \}$
are connected.

We also see by the way (this remark is not indispensable) 
that the inverse image ${\mathcal X}^0$
of 
$\{ b' = 0
\}$ has at most two connected  components, since it is the product of the  half-line $ \{ b > 0 \}$
with a central quadric $ \{ det B = a > 0\} \cap \{ b' = 0 \}$  in $\R^3$ with 
possible signatures $(+,+,-)$,$(+,- ,-)$,$(+,-)$.

In order to prove that ${\mathcal X}$ is connected, 
let us observe that the space  ${\mathcal X}$  we are considering is the
intersection of a  real quadric $Q$ in $\R^5$ with an open set.
Moreover, the quadric $Q$ is centred in the origin and the associated
quadratic form has negativity index, respectively positivity index,
at least two. Therefore $Q$ is everywhere of pure dimension $4$,
thus the  closed set ${\mathcal X}^0$ is in the closure
of the complement ${\mathcal X} \setminus {\mathcal X}^0$,
and it suffices to show that there is a point $p$ of ${\mathcal X}^0$
and a neighbourhood of $p$ meeting both ${\mathcal X}^+$
and ${\mathcal X}^-$.

But if this were not so, in the points of ${\mathcal X}^0$
the linear form $ b'$ restricted to ${\mathcal X}$
would vanish with its derivatives everywhere,
which implies that ${\mathcal X}^0$ would be a linear subspace
counted with multiplicity two, which is not the case.

We have thus shown that if $L=0$, ${\mathcal X}$ is connected. 

Assume now that $L \neq 0$. In this case we also have to consider equations
(\ref{Lanti}):    
 
$$L^+_+ B = -b L^-_-,$$
$$L^-_+ B = b^{-1} L^+_-.$$

If $L^+_+$, $L^-_+$ are linearly independent we can uniquely determine $B$ 
$$B = \left( \begin{array}{c} 
L^+_+\\ 
L^-_+
\end{array} \right)^{-1} \left( \begin{array}{c} 
-b L^-_-\\ 
b^{-1} L^+_-
\end{array} \right),$$
where $\left( \begin{array}{c} 
L^+_+\\ 
L^-_+
\end{array} \right)$ denotes the $(2 \times 2)$ matrix whose first row is
$L^+_+$ and whose second row is $L^-_+$ and $\left( \begin{array}{c} 
-b L^-_-\\ 
b^{-1} L^+_-
\end{array} \right)$ is the $(2 \times 2)$ matrix whose first row is
$-b L^-_-$ and whose second row is $b^{-1} L^+_-$.

Then $det(B) = - \frac{det (\left( \begin{array}{c} 
L^-_-\\ 
L^+_-
\end{array} \right) )} {det(\left( \begin{array}{c} 
L^+_+\\ 
L^-_+
\end{array} \right) )}= : \alpha$. 

So we must have $\alpha >0$, and equation (\ref{l=0}) is of
the form
$$a_- - \alpha a_+ = -b(c_1 b + c_2 b^{-1}),$$
where $c_1, c_2$ are functions depending only on $D$ and $L$. 

So we have to solve an equation of the form
$$c_1  b^2 + c =0,$$
and we  have at most one positive solution, thus ${\mathcal X}$ is connected.

Suppose now $L^-_+ = 0$. Then if there is a solution
also $L^+_- =0$ and we may assume without
loss of generality that $L^+_+ = (1, 0)$ (since $L \neq 0$). 

If we then set $B = (b_{ij})$ our equations (\ref{Lanti}) reduce to 
$$(b_{11}, b_{12}) =  -b L^-_-,$$
and equation (\ref{l=0}) becomes

$$a_- - a_+ det(B) = -b(c_0 b_{21} + c_2 b_{22} + b c_1),$$
where the $c_j$'s  only depend on $D$ and $L$. 

Moreover, we can write $det(B) = b \rho$, where $\rho = l_0 b_{21} + l_2 b_{22}$
 and the $l_j$'s  depend only on  $L$. 
So we have an equation in $ b, b_{21} , b_{22}$ 
$$a_- - a_+ b \rho = -b(c_0 b_{21} + c_2 b_{22} + b c_1),$$
which is an equation  of the form 
$$b^2 c_1 + b(a_0 b_{21} + a_2 b_{22}) + a_- = 0,$$
where the $a_j$'s  depend only on $L$ and $D$. 

The condition $det(B) > 0$ becomes $\rho >0$, which is a linear inequality in
$b_{21}, b_{22}$.  If we fix $b  \in  \R_+$, our solution is the intersection
of the  line in $\R^2$ given by 
$$a_0 b_{21} + a_2 b_{22}  + \frac{a_- + b^2 c_1}{b} =0$$
with the half plane $\rho >0$. 

To simplify things, we introduce a linear form $ \tau$ in the $\R^2$
with coordinates $b_{21}, b_{22}$, namely
$$ \tau : =  a_0 b_{21} + a_2 b_{22}.$$

There are three cases: 

1) the linear forms  $\rho  ,\tau$ are independent

2) $ \tau = c \rho$ for a constant $c \neq 0$.

3) $ \tau = 0$.

In case 1) we have  $ \tau = - \frac{a_- + b^2 c_1}{b}$, $ \rho > 0, b >0$,
and $\mathcal X$ is homeomorphic to a quadrant in $ \R^2$.

In case 2)  we have that $\mathcal X$ is a product of $\R$ with
the set $ \{ (\rho, b) |  \rho >0 , b > 0,   \rho = - \frac{a_- + b^2 c_1}{b c} \} \subset \R^2$,
which is diffeomorphic to the  interval
 $ \{ ( b) |   b > 0,   0 >  c (a_- + b^2 c_1 ) \} \subset \R$.

In case 3) $\mathcal X$ is a product of $\R$ with
the set $$ \{ (\rho, b) |  \rho >0 , b > 0,   0 = a_- + b^2 c_1 \} \subset \R^2.$$
which is a half-line in $\R^2$,

Therefore  $\mathcal X$ is connected in all three cases.

Assume now $L^-_+ \neq 0 $, $L^+_+ = \beta L^-_+$, $\beta \neq 0$. 
Then 
$$L^+_+ B = -b L^-_-,$$
$$L^-_+ B = b^{-1} L^+_-,$$

yield an equation of the form 
$$b^2 L^-_- = - \beta L^+_-,$$
and if we solve for $b$ we find at most one positive solution $\hat{b}$. 

Without loss of generality we can assume $L^+_+ = (1,0)$ and the equation 
$L^+_+ B = -\hat{b} L^-_- $
allows us to determine the first row of $B$. 

Hence both $det(B)$ and equation (\ref{l=0}) are 
polynomials of degree one  in $b_{21}$, $b_{22}$,
therefore 
the intersection of $\{ det(B)>0 \}$ with the set where the equation 
(\ref{l=0}) is satified is connected (a half-line) or
empty.

Finally,  if $L^+_+ =0$, then also $L^-_- =0$ and we may assume 
w.l.o.g. $L^-_+ = (1,
0)$.  
Therefore we have 
$$(b_{11}, b_{12}) = b^{-1} L^+_- =: b^{-1} (l_1, l_2),$$
and since then $a_- - det(B) a_+ = a_- - a_+ b^{-1}( l_1 b_{22} - l_2 b_{21}) $
our equation becomes
$$ a_- - a_+ b^{-1}( l_1 b_{22} - l_2 b_{21}) = 
 -b(b^{-1}(d_{11} l_2 - d_{21} l_1) + d_{12} b_{22} -
d_{22} b_{21}),$$
and the condition $det(B) >0$ can be written as $l_1 b_{22} - l_2 b_{21} >0$. 
 
So we get 

$$b_{22}(d_{12} b^2 - a_+ l_1) + b_{21}(- b^2 d_{22} + a_+ l_2) = 
-b(a_- - d_{21} l_1 + d_{11} l_2),$$
$$l_1 b_{22} - l_2 b_{21} >0.$$ 
 
We define functions 
$$ x (  b_{22}, b_{21}  ) : = l_1 b_{22} - l_2 b_{21} $$
$$ y ( b, b_{22}, b_{21}  ) : = b_{22}(d_{12} b^2 - a_+ l_1) + b_{21}(- b^2 d_{22} + a_+ l_2) .$$
These are linear functions of  $(b_{22}, b_{21})$
with determinant $  b^2 (d_{12} l_2 - d_{22} l_1)$.

Case 1) : $ (d_{12} l_2 - d_{22} l_1) \neq 0$.

Then $ ( b, b_{22}, b_{21}  ) \ra  (b, x (  b_{22}, b_{21}  ),  y ( b, b_{22}, b_{21}  ))$
is a self-diffeomorphism of $ \R_+ \times \R^2$, which transforms  
 ${\mathcal X}$ into the set  ${\mathcal X}' \subset \R^3$,
 ${\mathcal X}' = \{ (b,x,y) | b >0, x > 0, y = -b(a_- - d_{21} l_1 + d_{11} l_2) \}$,
which is diffeomorphic to the quadrant  $ \{ (b,x) | b >0, x > 0 \} $,
which is clearly connected.

Case 2) : $ (d_{12} l_2 - d_{22} l_1) = 0$.

In this case 
$$ y ( b, b_{22}, b_{21}  ) : = b_{22}(d_{12} b^2 - a_+ l_1) + b_{21}(- b^2 d_{22} + a_+
l_2)  =  ( c b^2 - a_+ ) \ x (  b_{22}, b_{21}  ) .$$  

There is a  linear self-diffeomorphism of $\R^3$ such that  
$$ ( b, b_{22}, b_{21}  ) \ra  (b, x ( b_{22}, b_{21}  ),  z ( b_{22}, b_{21}  ))).$$

It carries ${\mathcal X}$ into the set  

$ \{ (b,x,z) | b >0, x > 0, ( c b^2 - a_+ ) x = -b(a_- - d_{21} l_1 + d_{11} l_2) \}$,
which is the product of $\R$ with the subset of $\R^2$

$${\mathcal Y} : = \{ (b,x) | b >0, x > 0, ( c b^2 - a_+ ) x = -b(a_- - d_{21} l_1 + d_{11} l_2) \}.$$

Case 2.1) : $(a_- - d_{21} l_1 + d_{11} l_2) = 0$.

Then we get $ ( c b^2 - a_+ )= 0  $, which has at most one positive solution,
and ${\mathcal Y}$ is diffeomorphic to the half-line $ \{ x > 0 \}$.

Case 2.2) : $(a_- - d_{21} l_1 + d_{11} l_2) \neq 0$, $ c=0$.

In this subcase either  ${\mathcal Y}$  is the empty set, or 
 ${\mathcal Y}$  is diffeomorphic to $ \{ (b \in \R) | b >0 \}$.

Case 2.3) : $(a_- - d_{21} l_1 + d_{11} l_2) \neq 0$, $ c \neq 0$.

In this subcase we get that  ${\mathcal Y}$ is diffeomorphic to the set
 $ \{ b | b  > 0, \  c b^2  <  a_+\}$ or to the set $ \{ b | b  > 0, \  c b^2  >  a_+\}$.

We are done, since  in both cases either ${\mathcal Y}$ is empty or
it is diffeomorphic to an interval.

\qed

From the above theorem follows now easily

{\bf Theorem C }{\em   Assume again $ d=1, m=2$, that $A$ is nondegenerate and 
fix the orbifold fundamental group $\Pi_{\sigma}$.
Then the real structures for torus bundles in the Appell  Humbert family
are parametrized by a connected family.}

\Proof
It suffices to combine the connectedness result that we have just proven 
in Theorem 4.9 with theorem 4.3, with Remark \ref{top},
and  with the fact that ${\mathcal T}'{\mathcal B}_{A}^{\R}(A_1,A_2,L)$
is connected if and only if ${\mathcal T}{\mathcal B}_{A}^{\R}(A_1,A_2,L)$
is connected.

\qed

As an immediate corollary we have 

{\bf Theorem D }{\em   Same assumptions as in theorem B: $ d=1, m=2$, and 
 i), ii), iii) are satisfied.

Then simplicity holds for  real  torus bundles in the Appell  Humbert family.}

\vfill

\noindent
{\bf Author's address:}

\bigskip

\noindent
Prof. Fabrizio Catanese\\
Lehrstuhl Mathematik VIII\\
Universit\"at Bayreuth\\
     D-95440, BAYREUTH, Germany

e-mail: Fabrizio.Catanese@uni-bayreuth.de

\noindent
Dr. Paola Frediani\\
Dipartimento di Matematica\\
Universit\`a di Pavia\\
     I-27100 Pavia, Italy

e-mail: paola.frediani@unipv.it

\end{document}